\documentclass[12pt]{amsart}
\usepackage{fullpage,amsmath,xypic, amsthm,amsfonts,amssymb,graphicx,setspace,amscd,float,hyperref,enumitem,}
\usepackage[mathscr]{euscript}

\newtheorem{thm}{Theorem}[section]
\newtheorem{lem}[thm]{Lemma}
\newtheorem{qst}[thm]{Question}
\newtheorem{prop}[thm]{Proposition}

\theoremstyle{definition}
\newtheorem{df}[thm]{Definition}
\newtheorem{rk}[thm]{Remark}
\newtheorem{com}[thm]{Comment}

\newtheorem*{mainthm}{Theorem}

  \newcommand{\NN}{\mathbb{N}}

  \newcommand{\RR}{\mathbb{R}}

  \newcommand{\ZZ}{\mathbb{Z}}

\newcommand{\ind}{\mbox{ind}}

\newcommand{\gind}{\mbox{ind}_{\rm geom}}

\newcommand{\Ind}{\mbox{ind}}

\begin{document}

\title{Lone Axes in Outer Space}
\author{Lee Mosher and Catherine Pfaff}

\address{\tt Department of Mathematics, Rutgers Newark
  \newline \indent  {\url{http://www.andromeda.rutgers.edu/~mosher}}, } \email{\tt mosher@andromeda.rutgers.edu}
  
\address{\tt Department of Mathematics, University of California at Santa Barbara 
  \newline \indent  {\url{http://www.math.uni-bielefeld.de/~cpfaff/}}, } \email{\tt cpfaff@math.ucsb.edu}

\begin{abstract}
Handel and Mosher define the axis bundle for a fully irreducible outer automorphism in \cite{hm11}. In this paper we give a necessary and sufficient condition for the axis bundle to consist of a unique periodic fold line. As a consequence, we give a setting, and means for identifying in this setting, when two elements of an outer automorphism group $Out(F_r)$ have conjugate powers.
\end{abstract}

\thanks{The first author is supported by a grant from the National Science Foundation. The second author was supported first by the ARCHIMEDE Labex (ANR-11-LABX- 0033) and the A*MIDEX project (ANR-11-IDEX-0001-02) funded by the ``Investissements d'Avenir'' French government program managed by the ANR. She is secondly supported by the CRC701 grant of the DFG, supporting the projects B1 and C13 in Bielefeld.}  

\date{}
\maketitle

\section{Introduction}

We let $Out(F_r)$ denote the outer automorphism group for a rank $r$ free group \emph{$F_r$}. Culler and Vogtmann \cite{cv86} defined a topological space $CV_r$, \emph{Outer Space}, on which $Out(F_r)$ acts properly with finite stabilizers, in analogy with the action of each mapping class group on its Teichm\"{u}ller space (see \cite{flp79}). In fact, the action of each $Out(F_r)$ on its Outer Space $CV_r$ has indeed proved to possess many of the same characteristics as the action of a mapping class group on its Teichm\"{u}ller space. For example, Levitt and Lustig \cite{ll03} proved that, as with a pseudo-Anosov acting on Teichm\"{u}ller space, each ``fully irreducible'' $\varphi \in Out(F_r)$ acts with North-South dynamics on the natural compactification $\overline{CV_r}$ of $CV_r$. A fully irreducible outer automorphism is the most commonly used analogue to a pseudo-Anosov. An element $\varphi \in Out(F_r)$ is \emph{fully irreducible} if no positive power $\varphi^k$ fixes the conjugacy class of a proper free factor of $F_r$.

Recall that points of Outer Space can be described as marked metric graphs up to isometry, by which we mean graphs whose fundamental group has been identified with the free group in a basepoint-free manner and who have lengths assigned to their edges (generally assumed to sum to one). As in \cite{hm11}, one can call a point $\Gamma$ in Outer Space a \emph{train track} for $\varphi$ when there exists an affine train track representative $g \colon \Gamma \to \Gamma$. An affine train track representative is a train track representative, in the sense of \cite{bh92}, such that each open interval inside each edge is stretched by a constant factor equal to the dilitation of $\varphi$. In \cite{hm11} Handel and Mosher answered the question, posed by Vogtmann: ``Is the set of train tracks for an irreducible automorphism contractible?'' They do so by defining, for a nongeometric fully irreducible $\varphi \in Out(F_r)$, its axis bundle, which they also show is a closed subset $\mathcal{A}_{\varphi}$ in $CV_r$ proper homotopy equivalent to a line, invariant under $\varphi$, and such that the two ends of $\mathcal{A}_{\varphi}$ limit on the repeller and attractor of the source-sink action of $\varphi$ on $\overline{CV_r}$. Outer automorphisms induced by homeomorphisms of compact surfaces are called \emph{geometric}, and are usually primarily studied as surface homeomorphisms. If $\varphi$ is a nongeometric, fully irreducible outer automorphism, then~\cite[Theorem 1.1]{hm11} gives three equivalent definitions of the axis bundle (see Subsection \ref{SS:AB} below), the third of which is  $\mathcal{A}_{\varphi} = \overline{\bigcup_{k=1}^{\infty} TT(\varphi^k)}$, where $TT(\varphi^k)$ is just the set of train track graphs for $\varphi^k$.

Unlike in the situation of a loxodromic isometry acting on hyperbolic space or of a pseudo-Anosov mapping class acting on Teichm\"{u}ller space, it appears that there is in general no natural axis for a fully irreducible outer automorphism acting on Outer Space and that the axis bundle is a good natural analogue, in spite of in general being so far from a single axis as to actually be multi-dimensional. Handel-Mosher, via a list of questions in \cite{hm11}, and Bridson-Vogtmann \cite[Question 3]{bv06}, more directly, ask:

\begin{qst}
Describe the geometry of the axis bundle for a fully irreducible outer automorphism acting on Outer Space.
\end{qst}

What we accomplish in this paper is to determine when a fully irreducible outer automorphism behaves more like a pseudo-Anosov mapping class by having an axis bundle that is just a single axis. Not only does this give a partial solution to the conjugacy problem for outer automorphisms of free groups, but it allows one to ``read off'' from an axis all train track representatives for the automorphism. Section 4 is dedicated entirely to explaining several applications of our main theorem.

The condition we prove for a unique axis relies on the ideal Whitehead graph of \cite{hm11}. The condition also relies on the \cite{gjll} rotationless index $i(\varphi)$ for a fully irreducible $\varphi \in Out(F_r)$. One can think of the rotationless index as a sum of terms, each of which records the number of vertices in a component of the ideal Whitehead graph. As originally defined, the rotationless index also records the branching behavior of the attracting tree, $T^+_{\varphi}$, for the source-sink action of $\varphi$ on $\overline{CV_r}$. Unlike in the surface case where one has the Poincar\'e-Hopf index equality, Gaboriau, J\"ager, Levitt, and Lustig proved in \cite{gjll} that there is instead a rotationless index inequality $0 > i(\varphi) \geq 1-r$ that each fully irreducible $\varphi \in Out(F_r)$ satisfies. (Here we have rewritten the inequality using the \cite{p12d} rotationless index definition, revised to be invariant under taking powers and to have its sign be consistent with the mapping class group case.)

What we prove in Theorem \ref{t:uniqueaxis} of Section \ref{S:Proofs} is a necessary and sufficient condition for an ageometric fully irreducible outer automorphism to have a unique axis. One may note that examples of fully irreducibles satisfying the conditions of Theorem \ref{t:uniqueaxis} can be found in \cite{p12c} and \cite{p12d} and it was in fact proved later, in \cite{kp14}, that satisfying these conditions is generic along a particular ``train track directed'' random walk.

\begin{mainthm}\textbf{\ref{t:uniqueaxis}} \emph{The axis bundle of an ageometric, fully irreducible outer automorphism $\varphi \in Out(F_r)$ is a unique axis precisely if both of the following two conditions hold:
~\\
\vspace{-6mm}
\begin{enumerate}
\item the rotationless index satisfies $i(\varphi) = \frac{3}{2}-r$ and 
\item no component of the ideal Whitehead graph $\mathcal{IW}(\varphi)$ has a cut vertex.
\end{enumerate}}
\end{mainthm}

The rotationless index is always a negative half-integer. Thus, one may observe that $\frac{3}{2}-r$ is as close to equalling the bound of $1-r$ as possible without actually equalling it. As equality is achieved precisely in the case of geometric and parageometric outer automorphisms, this means that $\frac{3}{2}-r$ is the bounding ageometric rotationless index.

Given a nongeometric, fully irreducible $\varphi \in Out(F_r)$, we let $ST(\varphi)$ denote the set of train track graphs for $\varphi$ on which there exists a \emph{fully stable} train track representative for $\varphi$, meaning that each power is stable in the sense of \cite{bh92}. We then define the \emph{stable axis bundle} as $\mathcal{SA}_{\varphi} = \overline{\bigcup_{k=1}^{\infty} ST(\varphi^k)}$, see Section \ref{S:SAB}. The stable axis bundle was introduced in~\cite[Section 6.5]{hm11} as an object of interest. Our approach to proving Theorem \ref{t:uniqueaxis} involves a study of the stable axis bundle, as proposed in \cite{hm11}.

\begin{mainthm}\textbf{\ref{t:stablebundle}} \emph{Suppose $\varphi \in Out(F_r)$ is ageometric and fully irreducible. Then the stable axis bundle $\mathcal{SA}_{\varphi}$ is a unique axis if and only if the rotationless index satisfies $i(\varphi) = \frac{3}{2}-r$. In that case it is a unique periodic fold line.}
\end{mainthm}

The connection between the main theorem and its generalization is that, for an ageometric, fully irreducible $\varphi \in Out(F_r)$ with rotationless index $i(\varphi) = \frac{3}{2}-r$, the ideal Whitehead graph $\mathcal{IW}(\varphi)$ not having cut vertices is equivalent to the stable axis bundle in fact being the entire axis bundle. We exploit here constructions of~\cite[Lemma 3.1]{hm11} where cut vertices lead to periodic Nielsen paths in unstable representatives.

\bigskip

\subsection{Remarks and further questions}

\vskip1pt

The proof of the main theorem will exhibit a sufficient condition for the axis bundle to be of dimension two or higher, namely the existence of an affine train track representative $g \colon \Gamma \to \Gamma$ having two or more illegal turns. In fact the axis bundle has \emph{local} dimension two or higher at the point represented by $\Gamma$; see the proof of Lemma~\ref{l:distinct}. This condition motivates some follow-up problems regarding the behavior of higher dimensional axis bundles:
\begin{itemize}
\item Is there a formula for the local dimension of the axis bundle at a point represented by an affine train track representative, correlating that definition with the structure of the set of illegal turns in some manner? Is there such a formula for general points of the axis bundle?
\item Is the local dimension a constant function on the axis bundle?
\end{itemize}
Regarding these questions, the proof of the main theorem gives some hints. For instance, it indicates that the answer might be ``no''; see Comment \ref{c:c1} following the statement of Lemma~\ref{l:ec}.

\bigskip

\subsection*{Acknowledgements}
 
The second author would like to thank the first for his inspiration, support, and ideas, as well as for remaining one of her favorite people to discuss math with. She would also like to thank Thierry Coulbois and Jerome Los for illuminating conversations and their interest in her work. Finally, she would like to thank Yael Algom-Kfir, Mladen Bestvina, Michael Handel, Ilya Kapovich, and Martin Lustig for their support and thoughtful conversations about the geometry of Outer Space.

\section{Preliminary definitions and notation}{\label{S:Prelims}}

\subsection{Train track representatives}

\vskip1pt

\begin{df}[Marked graphs and train track representatives]
Let $R_r$ be the $r$-petaled rose, i.e. the graph with precisely $r$ edges and one vertex. Recall from \cite{bh92}, for example, that a \emph{marked graph} is a connected graph $\Gamma$, with no valence $1$ or $2$ vertices, together with an isomorphism $\pi_1(\Gamma) \cong F_r$ defined via a homotopy equivalence (called the \emph{marking}) $\rho \colon \Gamma \to R_r$. Marked graphs $\rho \colon \Gamma \to R_r$ and $\rho' \colon \Gamma' \to R_r$ are considered equivalent when there exists a homeomorphism $h \colon \Gamma \to \Gamma'$ such that $\rho' \circ h$ is homotopic to $\rho$. A homotopy equivalence $g \colon \Gamma \to \Gamma$ of a marked graph $\Gamma$ is a \emph{train track representative} for $\varphi \in Out(F_r)$ if it maps vertices to vertices, $\varphi = g_{\ast} \colon \pi_1(\Gamma) \to \pi_1(\Gamma)$, and $g^k \mid _{int(e)}$ is locally injective for each edge $e$ of $\Gamma$ and $k>0$. 
\end{df}


Many of the definitions and notation for discussing train track representatives were established in \cite{bh92} and \cite{bfh00}. We recall some here. 

\begin{df}[Turns and gates]
Let $g \colon \Gamma \to \Gamma$ be a train track representative of $\varphi \in Out(F_r)$. By a \emph{direction} at a vertex $v$ we will mean a germ of initial segments of directed edges emanating from $v$. The definition can be extended to an interior point $x$ of an edge $e$ by defining a direction at $x$ to be a germ of open segments of $e$ with $x$ as a boundary point. \emph{$Dg$} will denote the direction map induced by $g$. We call a direction $d$ \emph{periodic} if $Dg^k(d)=d$ for some $k>0$. We call an unordered pair of directions $\{d_i, d_j\}$, based at the same point, a \emph{turn}. It is an \emph{illegal turn} for $g$ if $Dg^k(d_i) = Dg^k(d_j)$ for some $k$ and a \emph{legal turn} otherwise. Considering the directions of an illegal turn equivalent, one can define an equivalence relation on the set of directions at a vertex. Each equivalence class is called a \emph{gate}.

\emph{Directions} and \emph{turns} at a point $v$ in a simplicial tree $T$ can be analogously defined, as can the direction map.
\end{df}

\vskip10pt

\subsection{Periodic Nielsen paths and (fully) stable representatives.}

\vskip1pt

Throughout this subsection, $\varphi \in Out(F_r)$ is fully irreducible and $g\colon\Gamma\to\Gamma$ is a train track representative of $\varphi$. (Hence, in particular, $g$ is expanding and irreducible.)

\begin{df}[Periodic Nielsen paths and principal points]{\label{d:PNPs}}
We call a locally injective path \emph{tight}.
Recall \cite{bh92}, a nontrivial tight path $\rho$ in $\Gamma$ is called a \emph{periodic Nielsen path (PNP)} for $g$ if, for some $k$, $g^k(\rho)\simeq\rho$ rel endpoints. It is called a \emph{Nielsen path (NP)} if the period is one and an \emph{indivisible Nielsen path (iNP)} if it further cannot be written as a concatenation $\rho = \rho_1\rho_2$, where $\rho_1$ and $\rho_2$ are also NPs for $g$.

As in \cite{hm11}, we call a periodic point $v \in \Gamma$ \emph{principal} that either has at least three periodic directions or is an endpoint of a periodic Nielsen path. 
\end{df}

\begin{df}[Rotationless]
A train track representative is called \emph{rotationless} if every principal point is fixed and every periodic direction at each principal point is fixed. Note that the rotationless property puts no restrictions on the preperiodic, nonperiodic vertices. In~\cite[Proposition 3.24]{fh11} it is shown that one can define a fully irreducible outer automorphism to be \emph{rotationless} if and only if one (hence all) of its train track representatives are rotationless.
\end{df}

We will use the following, which tells us that rotationless powers always exist: 

\smallskip

\begin{prop}~\cite[Corollary 4.43]{fh11}
For each $r \geq 2$, there exists an $R(r) \in \NN$ such that $\varphi^{R(r)}$ is rotationless for each $\varphi \in Out(F_r)$.
\end{prop}

\begin{df}[Stable train track representatives]
Let $\varphi$ be a fully irreducible outer automorphism. \cite{bh92} gives an algorithm for finding a representative with the minimal number of Nielsen paths, such a representative is called a \emph{stable} representative. As in \cite{hm11}, we call a stable representative $g$ of a rotationless power $\varphi^R$ of $\varphi$ \emph{fully stable}.
\end{df} 

\begin{rk}
It would not effect the definition of $\cup ST(\varphi^k)$ if we also called a representative fully stable whose rotationless powers are fully stable, but we will generally mean a rotationless representative when we use the term ``fully stable.'' 
\end{rk}

\vskip10pt

\subsection{Culler-Vogtmann Outer Space $CV_r$ and the attracting tree $T_+$}

\vskip1pt

\begin{df}[Outer space $CV_r$]
One can describe points in Culler-Vogtmann \emph{Outer Space} $CV_r$ as equivalence classes of marked metric graphs, where a \emph{metric} here means an assignment of positive lengths (summing to one) to the edges of the graph. At times one will in fact need a stronger version of a metric obtained by choosing for each edge $e$ of $\Gamma$, of length $l(e)$, a map $j_e \colon [0,l(e)] \to e$ that restricts to a homeomorphism on edge interiors. One can use these maps to define a path metric on $\Gamma$. 

The \emph{volume} of $\Gamma$ is defined as $vol(\Gamma) := \sum\limits_{e \in E(\Gamma)} l(e)$. Hence, in $CV_r$, the volume of each graph is one. However, there is an unprojectivized version of Outer Space, denote $\widehat{CV_r}$, where we no longer require $vol(\Gamma)=1$. 
\end{df} 

\begin{rk} Many definitions can be given equally for $CV_r$ and $\widehat{CV_r}$. Hence, we sometimes blur the distinction.
\end{rk}

Lifting to the universal cover, one obtains an alternative definition of $CV_r$, used to describe the compactification. Points of compactified Outer Space $\overline{CV_r} = CV_r \cup \partial CV_r$ can be described as equivalence classes of minimal, very small $F_r$-actions on $\mathbb{R}$-trees, known as ``$F_r$-trees.'' The equivalence relation is $F_r$-equivariant homothety. Under this description of $\overline{CV_r}$, points of $CV_r$ itself correspond to the simplicial $F_r$-trees $T$ on which $F_r$ acts freely; up to equivalence such trees correspond bijectively to marked graphs via the relation of universal covering.

There are multiple equivalent descriptions of the standard topology on $CV_r$. We describe it via its (ideal) simplicial structure. For each marked graph $\rho \colon \Gamma \to R_r$ with $N$ edges, the set of metrics on $\Gamma$ gives an $(N-1)$-dimensional open simplex in $CV_r$:

$$\{(l_1,l_2,\dots,l_N)\mid l_k >0, \sum l_k=1 \}.$$

Where they exist, open faces of a cell can be obtained by assigning length zero to a subset of the edges or equivalently by collapsing the forest in $\Gamma$ consisting of those edges of length zero. Faces are missing where assigning length zero to edges changes the homotopy type of the graph.

The group $Out(F_r)$ acts on $CV_r$ from the right, where each $\varphi \in Out(F_r)$ acts by precomposing the marking with an automorphism representing $\varphi$. Given a fully irreducible $\varphi \in Out(F_r)$, the repeller and attractor for the action on $\overline{CV_r}$ are elements of $\partial CV_r$, thus $F_r$-trees. We denote the attracting tree in $\partial CV_r$ by $T_+^{\varphi}$, or just $T_+$, and the repelling tree by $T_-^{\varphi}$, or just $T_-$. 

\begin{df}[Attracting tree $T_+^{\varphi}$]{\label{d:AttractingTree}}
We recall from \cite{gjll} a concrete construction of the attracting tree $T_+^{\varphi}$ for a fully irreducible $\varphi \in Out(F_r)$. Let $g\colon \Gamma \to \Gamma$ be a train track representative of $\varphi$ and $\tilde{\Gamma}$ the universal cover of $\Gamma$ equipped with a distance function $\tilde{d}$ lifted from $\Gamma$. As the fundamental group, $F_r$ acts by deck transformations, hence isometries, on $\tilde{\Gamma}$. A lift $\tilde{g}$ of $g$ is associated to a unique automorphism $\Phi$ representing $\varphi$. In particular, for each $w \in F_r$ and $x \in \tilde{\Gamma}$, we have $\Phi(w)\tilde{g}(x)=\tilde{g}(wx)$. One can define the pseudo-distance $d_{\infty}$ on $\tilde{\Gamma}$ by $\lim_{k \to +\infty} d_k$, where

$$d_k(x,y)=\frac{d(\tilde{g}^k(x),\tilde{g}^k(y))}{\lambda^k}$$ 

\noindent for each $x,y \in \tilde{\Gamma}$. Then $T_+$ is the $F_r$-tree defined by identifying each pairs of points $x,y \in \tilde{\Gamma}$ such that $d_{\infty}(x,y)=0$.
\end{df}

\vskip10pt

\subsection{The attracting lamination $\Lambda_{\varphi}$ for a fully irreducible outer automorphism.}

\vskip1pt
While one can define the set of attracting laminations for any element of $Out(F_r)$ (see \cite{bfh97}), we give here only a definition yielding the unique (see~\cite[Lemma 1.12]{bfh00}) attracting lamination for a fully irreducible.

Let $\Gamma$ be a marked graph with universal cover $\tilde{\Gamma}$ and projection map $p \colon \tilde{\Gamma} \to \Gamma$. By a \emph{line} in $\tilde{\Gamma}$ we mean the image of a proper embedding of the real line $\tilde{\lambda} \colon \mathbb{R} \to \tilde{\Gamma}$. We denote by $\tilde{\mathcal{B}}(\Gamma)$ the space of lines in $\tilde{\Gamma}$ with the compact-open topology (generated by the open sets $\tilde{\mathcal{U}}(\tilde{\gamma}):=\{L \in \tilde{\mathcal{B}}(\Gamma) \mid \tilde{\gamma}$ is a finite subpath of $L\}$). A \emph{line} in $\Gamma$ is then the image of a projection $p \circ \tilde{\lambda} \colon \mathbb{R} \to \Gamma$ of a line $\tilde{\lambda}$ in $\tilde{\Gamma}$, where two lines are considered equivalent when they differ via precomposition by a homeomorphism of $\RR$. We denote by $\mathcal{B}(\Gamma)$ the space of lines in $\Gamma$ with the quotient topology induced by the natural projection map from $\tilde{\mathcal{B}}(\Gamma)$ to $\mathcal{B}(\Gamma)$. One can then define the $\mathcal{U}(\gamma):=\{L \in \mathcal{B}(\Gamma) \mid \gamma$ is a finite subpath of $L\}$, which generate the topology on $\mathcal{B}$. For a marked graph $\Gamma$, we say a line $\gamma$ in $\Gamma$ is \emph{birecurrent} if every finite subpath of $\gamma$ occurs infinitely often as an unoriented subpath in each end of $\gamma$.

\begin{df}[Attracting lamination $\Lambda_{\varphi}$]
Fix a fully irreducible $\phi \in Out(F_r)$ and consider any train track representative $g \colon \Gamma \to \Gamma$ for a fully irreducible $\varphi \in Out(F_r)$. Given any edge $e$ in $\Gamma$, there exists a $k>0$ such that the following is a sequence of nested open sets:
    $$\mathcal{U}(e) \supset \mathcal{U}(g^k(e)) \supset \mathcal{U}(g^{2k}(e)) \dots$$
The \emph{attracting lamination} $\Lambda_{\varphi}$ (or just $\Lambda$) for $\varphi$ is the set of birecurrent lines in the intersection. We often use the same notation for the total lift $\widetilde{\Lambda}$ of $\Lambda$ to the universal cover. The meaning should be clear from context.
\end{df}

This definition of $\Lambda$ is well-defined independent of the choice of train track representative; see \cite[Lemma 1.12]{bfh97} for proof. Notice that, once a basepoint lift is chosen in $\tilde{\Gamma}$, one can identify $\partial \tilde{\Gamma}$ with the hyperbolic boundary $\partial F_r$ of the free group. This allows one to identify $\tilde{\Lambda}$ with a set of unordered pairs of points in $\partial F_r$, by lifting $\Lambda$ via the projection $\partial F_r = \partial\tilde{\Gamma} \to \mathcal{B}(\Gamma)$. It follows that $\tilde{\Lambda}$ is also well-defined.

We may also define the realization of $\Lambda$ in a general point of Outer Space represented by a marked graph $\Gamma'$ with universal cover $\tilde{\Gamma'}$ and with a chosen basepoint in $\tilde{\Gamma'}$. Using the identifications $\partial \tilde{\Gamma} \approx \partial F_r \approx \partial \tilde{\Gamma'}$, we obtain an identification $\tilde{\mathcal{B}}(\Gamma) \approx \tilde{\mathcal{B}}(\Gamma')$, which identifies $\tilde{\Lambda} \subset \tilde{\mathcal{B}}(\Gamma)$ with a subset of $\tilde{\mathcal{B}}(\Gamma')$ which is the realization of $\tilde \Lambda$ in $\tilde \Gamma'$. Following by the projection $\tilde{\mathcal{B}}(\Gamma') \to \mathcal{B}(\Gamma')$, we obtain the realization of $\tilde{\Lambda}$ in $\Gamma'$.

\vskip10pt

\subsection{$\Lambda$-isometries and weak train tracks.}

\vskip1pt

Recall that a fully irreducible $\varphi \in Out(F_r)$ is \emph{geometric} if it is represented by a homeomorphism $f \colon M \to M$ of a compact surface with nonempty boundary, meaning that there exists a homotopy equivalence $h \colon R_r \to M$ with homotopy inverse $\bar h \colon M \to R_r$, such that the homotopy equivalence $\bar h f h$ is homotopic to a train track representative of $\varphi$. Let $\varphi \in Out(F_r)$ be a nongeometric fully irreducible with attracting lamination $\Lambda$.

\begin{df}[$\Lambda$-isometry]
For a free, simplicial $F_r$-tree $T$, a \emph{$\Lambda$-isometry} on $T$ is an $F_r$-equivariant map $f_T \colon T \to T_+$ such that, for each leaf $L$ of $\Lambda$ realized in $T$, the restriction of $f_T$ to $L$ is an isometry onto a bi-infinite geodesic in $T_+$.
\end{df} 

\begin{df}[Weak train track]
A \emph{weak train track (WTT)} for $\varphi$ is a free, simplicial $F_r$-tree $T$ on which a $\Lambda$-isometry exists.
\end{df} 

\vskip10pt

\subsection{The axis bundle.}\label{SS:AB}

\vskip1pt

Three equivalent definitions of the axis bundle $\mathcal{A}_{\varphi}$ for a nongeometric fully irreducible outer automorphism $\varphi \in Out(F_r)$ are given in \cite{hm11}. We will use all three definitions here and thus remind the reader in this subsection of each of them. We say a few words with regard to their equivalence in Subsection \ref{SS:ABEquiv}.

\begin{df}[Fold lines]
A \emph{fold line} in $CV_r$ is a continuous, injective, proper function $\mathbb{R} \to CV_r$ defined by \newline
\noindent 1. a continuous 1-parameter family of marked graphs $t \to \Gamma_t$ and \newline
\noindent 2. a family of homotopy equivalences $h_{ts} \colon \Gamma_s \to \Gamma_t$ defined for $s \leq t \in \mathbb{R}$, each marking-preserving, \newline
\indent satisfying:
~\\
\vspace{-\baselineskip}
\begin{description}
\item [\emph{Train track property}] $h_{ts}$ is a local isometry on each edge for all $s \leq t \in \mathbb{R}$. 
\item [\emph{Semiflow property}] $h_{ut} \circ h_{ts} = h_{us}$ for all $s \leq t \leq u \in \mathbb{R}$ and $h_{ss} \colon \Gamma_s \to \Gamma_s$ is the identity for all $s \in \mathbb{R}$.
\end{description}
\end{df}

\noindent \textbf{Axis Bundle Definition I.}  $\mathcal{A}_{\varphi}$ is the union of the images of all fold lines $\mathcal{F} \colon \mathbb{R} \to CV_r$ such that $\mathcal{F}$(t) converges in $\overline{CV_r}$ to $T_{-}^{\varphi}$ as $t \to -\infty$ and to $T_{+}^{\varphi}$ as $t \to +\infty$.

\bigskip

\noindent \textbf{Axis Bundle Definition II.} $\mathcal{A}_{\varphi}$ is the set of weak train tracks for $\varphi$, i.e. 

$\mathcal{A}_{\varphi} =$ \{free simplicial $F_r$-trees $T \in CV_r \mid$ there exists a $\Lambda$-isometry $f_T \colon T \to T_{+}$\}.

\bigskip

\noindent \textbf{Axis Bundle Definition III.} $\mathcal{A}_{\varphi} = \overline{\bigcup_{k=1}^{\infty} TT(\varphi^k)}$, where $TT(\varphi^k)$ is the set of train track graphs for $\varphi^k$ and the closure is taken in $CV_r$.

\begin{df} We denote by $\widehat{\mathcal{A}_{\varphi}}$ the full lift of $\mathcal{A}_{\varphi}$ to unprojectivized Outer Space $\widehat{CV_r}$.
\end{df}

\vskip10pt

Several crucial properties of the axis bundle are recorded in \cite[Theorem 6.1, Lemma 6.2]{hm11}. We summarize a few here as Proposition \ref{P:6.2}. Given a point $T \in CV_r$, $Len(T) :=vol(T/F_r)$.

\begin{prop}[\cite{hm11}]\label{P:6.2} Let $\varphi \in Out(F_r)$ be a nongeometric fully irreducible. Then the map $Len \colon \widehat{\mathcal{A}_{\varphi}} \to (0, \infty)$ is a surjective and $\varphi$-equivariant homotopy equivalence where $\varphi$ acts on $(0, +\infty)$ by multiplication by $\frac{1}{\lambda}$.
\end{prop}

\vskip10pt

\subsection{A bit on equivalence of the axis bundle definitions.}\label{SS:ABEquiv}

\vskip1pt

The equivalence of the three axis bundle definitions is proved in~\cite[Theorem 1.1]{hm11}. We explain here briefly connections frequently used. In particular, we show one obtains from a train track representative $g \colon \Gamma \to \Gamma$ of a nongeometric fully irreducible $\varphi \in Out(F_r)$ both a $\Lambda$-isometry $g_{\infty} \colon \tilde{\Gamma} \to T_+$ and a ``periodic'' fold line.

\subsubsection{$\Lambda$-isometries from train track maps.}\label{SSS:DirectLimitMap}

\begin{df}[$g_{\infty}$]
Let $g \colon \Gamma \to \Gamma$ be a train track representative of a nongeometric fully irreducible $\varphi \in Out(F_r)$. We return to the construction of Definition \ref{d:AttractingTree}. We let $T_k$ denote the simplicial $F_r$-tree obtained from $\tilde{\Gamma}$ by identifying each $x,y \in \tilde{\Gamma}$ with $d_k(x,y)=0$ and then equipping the quotient graph with the metric induced by $d_k$. Then, for each $i$, a basepoint-preserving lift of $g$ induces a basepoint-preserving $F_r$-equivariant map $\tilde{g}_{i+1,i} \colon T_i \to T_{i+1}$ restricting to an isometry on each edge. We obtain a direct system $\tilde{g}_{j,i} \colon T_i \to T_j$ defined inductively by $\tilde{g}_{j,i}=\tilde{g}_{j,j-1} \circ \tilde{g}_{j-1,i}$. Then the $\Lambda$-isometry $g_{\infty} \colon \tilde{\Gamma} \to T_+$ is the direct limit map.

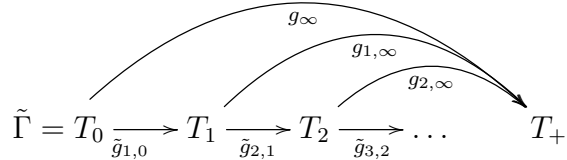
\begin{figure}[here]
\[
\xymatrix{\tilde{\Gamma}=T_0 \ar[r]_{\tilde{g}_{1,0}} \ar@/^4pc/[rrrr]_{g_{\infty}} & T_1 \ar[r]_{\tilde{g}_{2,1}} \ar@/^3pc/[rrr]_{g_{1,\infty}} & T_2 \ar[r]_{\tilde{g}_{3,2}} \ar@/^2pc/[rr]_{g_{2,\infty}} & \dots   & T_+ \\}
\]
\end{figure}
\end{df}

We will use the following, which is \cite[Corollary 2.14]{hm11}.

\begin{prop}[\cite{hm11}]\label{P:C2.14}
Let $g \colon \Gamma \to \Gamma$ be a train track representative of a nongeometric fully irreducible $\varphi \in Out(F_r)$, let $\tilde{\Gamma}$ be the universal cover, and let $[\tilde{x},\tilde{y}]$ be the tight path from $\tilde{x}$ to $\tilde{y}$. Then $g_{\infty} \colon \tilde{\Gamma} \to T_+$ is a surjective equivariant map such that, for all $\tilde{x},\tilde{y} \in \tilde{\Gamma}$, the following are equivalent:
\begin{enumerate} 
\item $g_{\infty}(\tilde{x})=g_{\infty}(\tilde{y})$
\item there exist $k \geq 0$ such that $g^k_{\#}([\tilde{x},\tilde{y}])$ is either trivial or a Nielsen path.
\end{enumerate}
In particular, $g_{\infty}$ restricts to an isometry on all legal paths.
\end{prop}

\begin{rk}[Realizing lamination leaves in $T_+$] Since lamination leaves are legal, Proposition \ref{P:C2.14} allows one to describe how lamination leaves are realized in $T_+$.
\end{rk}

\vskip10pt

\subsubsection{Periodic fold lines.}\label{SSS:PeriodicFoldLines}

\vskip1pt

Stallings introduced ``folds'' in \cite{s83}. Let $g \colon \Gamma \to \Gamma'$ be a homotopy equivalence of marked graphs. Let $e_1' \subset e_1$ and $e_2' \subset e_2$ be maximal, initial, nontrivial subsegments of edges $e_1$ and $e_2$ emanating from a common vertex and satisfying that $g(e_1')=g(e_2')$ as edge paths and that the terminal endpoints of $e_1'$ and $e_2'$ are in $g^{-1}(\mathcal{V}(\Gamma))$. Redefine $\Gamma$ to have vertices at the endpoints of $e_1'$ and $e_2'$ if necessary. One can obtain a graph $\Gamma_1$ by identifying the points of $e_1'$ and $e_2'$ that have the same image under $g$, a process we will call \emph{folding}. Stallings \cite{s83} also showed that if $g \colon \Gamma \to \Gamma'$ is tight, then $g$ factors as a composition of folds and a final homeomorphism. We call such a decomposition a \emph{Stallings fold decomposition}. It can be obtained as follows: At an illegal turn for $g\colon \Gamma  \to \Gamma'$, one can fold two maximal initial segments having the same image in $\Gamma'$ to obtain a map $\mathfrak{g}_1 \colon \Gamma_1 \to \Gamma'$ of the quotient graph $\Gamma_1$. The process can be repeated for $\mathfrak{g}_1$ and recursively. If some $\mathfrak{g}_k \colon \Gamma_{k-1} \to \Gamma$ has no illegal turn, then $\mathfrak{g}_k$ will be a homeomorphism and the fold sequence is complete.

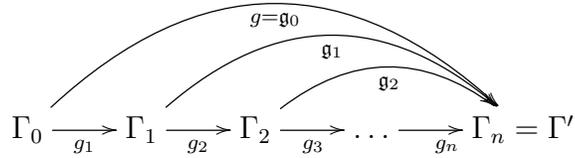
\begin{figure}[here]
\[
\xymatrix{\Gamma_0 \ar[r]_{g_1} \ar@/^4pc/[rrrr]_{g=\mathfrak{g}_0} & \Gamma_1 \ar[r]_{g_2} \ar@/^3pc/[rrr]_{\mathfrak{g}_1} & \Gamma_2 \ar[r]_{g_3} \ar@/^2pc/[rr]_{\mathfrak{g}_2} & \dots \ar[r]_{g_n}  & \Gamma_n=\Gamma' \\}
\]
\caption{Constructing a Stallings folds decomposition}\label{Fi:Sta}
\end{figure}

Notice that choices of illegal turns are made in this process and that different choices lead to different Stallings fold decompositions of the same homotopy equivalence.

When $\Gamma$ is a marked metric graph (of volume 1), we obtain an induced metric on each $\Gamma_k$, which we may renormalize to be again of volume 1.

In \cite{s89}, Skora interpreted a Stallings fold decomposition for a tight homotopy equivalence $g\colon \Gamma \to \Gamma'$ as a sequence of folds performed continuously. Repeating a Stallings fold decomposition of a train track representative for a fully irreducible defines a periodic fold line in Outer Space. The discretization of this fold line is depicted in Equation \ref{E:PeriodicFoldLines} below, where it should be noted that $\Gamma_{nk}=\frac{1}{\lambda^n}\Gamma_0 \cdot \varphi^n$, for each integer $n$.

\begin{equation}\label{E:PeriodicFoldLines}
\dots \xrightarrow{} \Gamma_0 \xrightarrow{g_1} \Gamma_1 \xrightarrow{g_2} \cdots \xrightarrow{g_K} \Gamma_K \xrightarrow{g_{K+1}} \Gamma_{K+1} \xrightarrow{g_{K+2}} \cdots \xrightarrow{g_{2K}} \Gamma_{2K} \xrightarrow{g_{2K+1}} \dots
\end{equation}

\vskip10pt

\subsection{Ideal Whitehead graphs and the rotationless index.}

\vskip1pt

We first explain for the reader more familiar with surface theory the ideal Whitehead graph and index list for a pseudo-Anosov $\psi$ on a closed surface $S$.  Suppose that $v$ is a $k$-pronged foliation singularity and $\tilde{\psi}$ is a lift of $\psi$ to the universal cover $\tilde{S}$ fixing a lift $\tilde{v}$ of $v$. Then $\tilde{v}$ in fact lies inside of a principal region $P$ for the invariant lamination. The index list entry for $v$ would be $1-\frac{k}{2}$ and the ideal Whitehead graph componenet would be a $k$-gon. Equivalently, the ideal Whitehead graph component is the polygon formed by the lamination leaf lifts bounding the principal region $P$.

We remind the reader of an \cite{hm11} definition for the ideal Whitehead graph of a nongeometric fully irreducible. One can reference \cite{p12a} and \cite{hm11} for alternative definitions of the ideal Whitehead graph and its outer automorphism invariance.

\begin{df}[Ideal Whitehead graph $\mathcal{IW}(\varphi)$]\label{d:IWG}
Let $\varphi \in Out(F_r)$ be a nongeometric fully irreducible with lifted attracting lamination $\tilde\Lambda$ (realized in $T_+$). To define the ideal Whitehead graph, start with the graph having a vertex for each distinct leaf endpoint and an edge connecting the vertices for the endpoints of each leaf. $\widetilde{\mathcal{IW}(\varphi)}$ is the union of the components with at least three vertices. $F_r$ acts freely, properly discontinuously, and cocompactly in such a way that the restriction to each component of $\widetilde{\mathcal{IW}(\varphi)}$ has trivial stabilizer. The \emph{ideal Whitehead graph} $\mathcal{IW}(\varphi)$ is the quotient under this action, which one can note has only finitely many components.
\end{df}

The ideal Whitehead graph has another interpretation in terms of singular leaves of the attracting lamination. For a fully irreducible $\varphi$, a leaf of the attracting lamination is \emph{singular} that shares a half-leaf with another leaf. Two half-leaves are \emph{asymptotic} if they share a common ray. $\widetilde{\mathcal{IW}(\varphi)}$ has a vertex for each asymptotic class of half-leaves of singular leaves and an edge for each singular leaf of $\tilde\Lambda$. The edges for two singular half-leaves share a vertex precisely when those half-leaves share an endpoint in the boundary. This definition also allows one to view $\widetilde{\mathcal{IW}(\varphi)}$ in any $F_r$-tree $T \in TT(\varphi^k)$. Many details of the correspondence of these views will be explained in Remark \ref{r:MixViews}.

\begin{rk}
As mentioned before, cut vertices of an ideal Whitehead graph yield representatives with PNPs. One can obtain such a representative from a stable representative by splitting open at cut vertices of the ideal Whitehead graph, as in~\cite[Lemma 4.3]{hm11}. We use this in particular in Lemma \ref{l:sisa}.
\end{rk}

The notion of an index $\ind{(\varphi)}$ for a $\varphi \in Out(F_r)$ was first introduced in \cite{gjll}. This notion is not in general invariant under taking powers. \cite{hm11} introduces the notion of a rotationless index (there just called the index sum) $i(\varphi)$ for a fully irreducible $\varphi \in Out(F_r)$. It follows from~\cite[Lemma 3.4]{hm11} that for rotationless nongeometric fully irreducible $\varphi \in Out(F_r)$, the two notions differ only by a change of sign.

\begin{df}[Index list and index sum]\label{d:IndexLists}
Let $\varphi \in Out(F_r)$ be a nongeometric fully irreducible and let $C_1, \dots, C_l$ be the connected components of $\mathcal{IW}(\varphi)$. For each $j$, let $k_j$ denote the number of vertices of $C_j$. The \emph{index list} for $\varphi$ is defined as

\begin{equation}\label{E:IndexLists} 
(i_1, \dots, i_j, \dots, i_l) = (1-\frac{k_1}{2}, \dots, 1-\frac{k_j}{2}, \dots, 1-\frac{k_l}{2}),
\end{equation} 

\noindent where the list is rewritten to be in increasing order of absolute values with repetitions allowed. The \emph{rotationless index} is then $i(\varphi) = \sum\limits_{j=1}^l i_j$.
\end{df}

One can obtain the index list (hence rotationless index) from any PNP-free rotationless train track representative $g \colon \Gamma \to \Gamma$. The $k_i$ in (\ref{E:IndexLists}) are replaced by the number of gates $k_i$ at the principal vertices $v_i \in \Gamma$. Since g is PNP-free, the principal vertices are precisely those periodic vertices with at least three gates. The index sum is therefore:

\begin{equation}\label{E:IndexSum}  
i(\varphi) = \sum\limits_{\text{principal vertices v}}(1-\frac{\#(\text{gates at v})}{2}).
\end{equation}

\vskip10pt

\subsection{Ageometrics.}

\vskip1pt

The division of the set of nongeometric fully irreducibles into ``ageometric'' and ``parageometric'' outer automorphisms could be considered to have evolved out of a series of papers. In \cite{gl95} Gaboriau and Levitt define the geometric index $\gind(T)$ for an $\RR$-tree equipped with a minimal, small $F_r$-action. They prove that the index satisfies the inequality $\frac{1}{2} \leq \gind(T) \leq r-1$, with the equality $\gind(T)=r-1$ realized precisely by ``geometric trees.'' In \cite{gjll} it is proved that, after replacing a fully irreducible $\varphi$ by a suitable positive power, one has $\gind(T_+^{\varphi}) = 2\Ind(\varphi)$. While it had been previously known that geometric fully irreducibles have geometric attracting tree, Levitt proved in \cite{l93} that even some nongeometric fully irreducibles have geometric attracting tree, hence creating a natural division of nongeometric fully irreducible outer automorphisms by their rotationless index. It will still be important for us that an ageometric fully irreducible $\varphi \in Out(F_r)$ can be characterized by satisfying $0 > i(\varphi) > 1-r$. However, we also give an equivalent definition in terms of PNPs. The equivalence follows from the fact proved in ~\cite[Theorem 3.2]{bf94} that for a fully irreducible $\varphi \in Out(F_r)$, the attracting tree $T_+^{\varphi}$ is geometric if and only if the ``stable'' train track representative of $\varphi$ contains a PNP. For interest's sake, we make one final remark that independently Handel and Mosher \cite{hm07} and Guirardel \cite{g05} gave a further characterization that a fully irreducible is geometric if and only if both the attracting tree and repelling tree are geometric.

\begin{df}[Ageometric]\label{D:Ageometric} 
A fully irreducible outer automorphism is \emph{ageometric} when a fully stable representative of a rotationless power has no Nielsen paths (closed or otherwise). 
\end{df}

\begin{rk}
Note that by ~\cite[Lemma 3.28]{fh11} every PNP of a rotationless train track representative is in fact an NP. Hence, a fully irreducible is ageometric if and only if one (hence every) fully stable representative of a (hence every) rotationless power has no NPs. 
\end{rk}

\subsection{Local decomposition of ideal Whitehead graphs.}

\vskip1pt

The following definitions are from \cite{hm11} and, as in \cite{hm11}, we assume throughout that $\varphi \in Out(F_r)$ is nongeometric fully irreducible, and $T_+ = T_+^{\phi}$, and $\Lambda = \Lambda_{\phi}$, and $f \colon T \to T_+$ is a $\Lambda$-isometry. We further let $\Gamma := T / F_r$.

There is a partial ordering on the axis bundle which relies on the splitting of $\mathcal{IW}(\varphi)$ into ``stable Whitehead graphs.'' 

We remind the reader of terminology and discussion found in~\cite[Section 3.3]{hm11}.

\begin{df}[Principal points in trees]\label{d:Principal}
Given a branch point $b$ of $T_+$, the lifted ideal Whitehead graph $\widetilde{\mathcal{IW}(\varphi)}$ has one component, which we denote $\mathcal{C}_b$, whose edges, realized as lines in $T_+$, all contain $b$. This relationship gives a one-to-one correspondence between components of $\widetilde{\mathcal{IW}(\varphi)}$ and branch points of $T_+$. Given a branch point $b$ of $T_+$, we let $\mathcal{C}_{b,T}$ denote the realization of $\mathcal{C}_b$ in $T$. This makes sense by viewing the ideal Whitehead graph in terms of the lamination leaves, as in Definition \ref{d:IWG}. We call a point $v$ in $T$ \emph{principal} for $f$ if there exists a branch point $b$ of $T_+$ such that $f(v)=b$ and $v$ is in $\mathcal{C}_{b,T}$.
\end{df} 

\begin{rk}
Notice that by definition, and by $f \colon T \to T_+$ being a $\Lambda$-isometry, that $f$, restricted to the principal points, is surjective onto the set of principal points, i.e. branch points of $T_+$. Thus, $f$ injects principal vertices to branch points if and only if it bijects principal vertices to branch points.
\end{rk}

\begin{df}[Basepoints and singular rays]
Let $l_T$ denote the realization in $T$ of a singular leaf $l$ of some $\mathcal{C}_b$. Then there exists a unique principal point $v$ of $T$, called the \emph{basepoint} of $l_T$, lying on $l_T$ and such that $f(v)=b$. Then $v$ divides $l$ into two rays. Rays obtained as such are called \emph{singular} rays.
\end{df} 

\begin{rk}{\label{r:MixViews}}
There exists a one-to-one correspondence between the set of singular rays in $T$ and the set of fixed directions based at principal vertices. The correspondence can be seen as follows (further details can be found in \cite{p12a}). Given a direction $d$ of an edge $E$ at a principal vertex $v$, the ray determined by $d$ is defined as $\tilde{R} = \cup_{j=0}^{j=\infty}\tilde{g}^j(E)$, where $\tilde{g}$ is a lift, fixing $v$, of a rotationless train track representative $g \colon \Gamma \to \Gamma$ for some $\varphi^k$.

Recall from Definition \ref{d:PNPs} that a principal point downstairs either has 3 periodic directions or is the endpoint of a PNP.
For a principal point downstairs, having at least three periodic directions, any lift is a principal point upstairs (with at least three singular rays determined by the singular directions, as in the previous paragraph). Additionally, since the directions are fixed, they give three distinct edges at $f(v)$, so that $f(v)$ is indeed a branchpoint of $T_+$.

Suppose that instead $v$, $w$ are the endpoints of an iNP $\rho$ for $g$. Then there exist at least two fixed directions at $v$ (one, which we call $E_1$, is a terminal edge of $\rho$) and at least two fixed directions at $w$ (one, which we call $E_2$, is the other terminal edge of $\rho$). For any lift $\tilde{\rho}$ of $\rho$ with terminal vertices $\tilde{v}$ and $\tilde{w}$, lifts of $v$ and $w$ respectively, $f$ sends $\tilde{v}$ and $\tilde{w}$ to a common point $b$ of $T_+$, which has at least three directions (one of which arises from the identification of $E_1$ and $E_2$ and the other two of which come from the distinct fixed directions at $v$ and $w$). Hence, $b$ is also a branchpoint $T_+$. Notice that the rays constructed above from $E_1$ and $E_2$ are also asymptotic in $T$.

All principal points and principal directions in $T$ arise in one of the two ways described above.
\end{rk}

\begin{df}[Stable Whitehead graphs $\mathcal{SW}(\tilde{v};T)$ and $\mathcal{SW}(v;\Gamma)$, and local Whitehead graph $\mathcal{LW}(v;\Gamma)$]{\label{d:WhiteheadGraphs}}
The \emph{local Whitehead graph} $\mathcal{LW}(\tilde{v};T)$ will have a vertex for each direction at $\tilde{v}$ and an edge connecting the vertices corresponding to the pair of directions $\{d_1,d_2\}$ if the turn $\{d_1,d_2\}$ is taken by the realization in $T$ of a leaf of $\Lambda_{\varphi}$. The \emph{stable Whitehead graph} $\mathcal{SW}(\tilde{v};T)$ at a principal point $\tilde{v}$ will be the subgraph of $\mathcal{LW}(\tilde{v};T)$ obtained by restricting to the periodic directions. Equivalently, the \emph{stable Whitehead graph} $\mathcal{SW}(\tilde{v};T)$ at a principal point $v$ can be identified with the graph having one vertex for each singular ray $\tilde{R}$ based at $v$ and an edge connecting the vertices corresponding to a pair of singular rays $\tilde{R_1}$, $\tilde{R_2}$ at $v$ if and only if $\tilde{R_1} \cup \tilde{R_2}$ is a singular leaf at $v$. One can reference \cite{hm11} for further details.

The \emph{local Whitehead graph} $\mathcal{LW}(v;\Gamma)$, at a point $v \in \Gamma$ again has a vertex for each direction at $v$ and an edge connecting the vertices corresponding to the pair of directions $\{d_1,d_2\}$ if the turn $\{d_1,d_2\}$ is taken by the realization in $\Gamma$ of a leaf of $\Lambda_{\varphi}$. And the \emph{stable Whitehead graph} $\mathcal{SW}(v;\Gamma)$ at a principal point $v$ will be the subgraph of $\mathcal{LW}(v;\Gamma)$ obtained by restricting to the periodic directions. Since each gate at a vertex contains precisely one periodic direction, one can equally give this definition in terms of gates.

Remark \ref{r:MixViews} explains the relationship between the different stable Whitehead graphs.

Notice that, given a train track map $g \colon \Gamma \to \Gamma$, the direction map $Dg$ induces a simplicial map on both the local Whitehead graph and the stable Whitehead graph. We again denote these maps by $Dg$.
\end{df}

\begin{rk}{\label{r:SWGinIWG}}
Each stable Whitehead graph $\mathcal{SW}(\tilde{v};T)$ sits inside of $\mathcal{IW}(\varphi)$: A vertex of $\mathcal{SW}(\tilde{v};T)$ corresponds to a singular leaf $\tilde{R}$ at $\tilde{v}$ and the endpoint of this ray corresponds to a vertex of $\mathcal{IW}(\varphi)$. An edge of $\mathcal{SW}(\tilde{v};T)$ corresponds to a singular leaf based at $\tilde{v}$. This leaf also gives an edge of $\mathcal{IW}(\varphi)$. 
\end{rk}

\begin{df}[Local decompositions and splitting]
In light of Remark \ref{r:SWGinIWG} (and \cite[Lemma 5.2]{hm11}), the ideal Whitehead graph, realized in $T$, can be written as the union of the stable Whitehead graphs at the principal points. We call this its \emph{local decomposition}.

Let $T,T'$ be weak train tracks with $\Lambda$-isometries $f \colon T \to T_+$ and $f' \colon T' \to T_+$. As in \cite{hm11}, one says $f$ \emph{splits as much as} $f'$ if the local decomposition $\mathcal{IW}(\varphi) = \bigcup \mathcal{SW}(v_j;T)$ is at least as fine as the local decomposition $\mathcal{IW}(\varphi) = \bigcup \mathcal{SW}(w_i;T')$. That is, for each principal vertex $v_j$ of $T$, there exists a principal vertex $w_i$ of $T'$ such that $\mathcal{SW}(v_j;T) \subset \mathcal{SW}(w_i;T)$, where the inclusion takes place in $\mathcal{IW}(\varphi)$, realized as a decomposition, as above. 
\end{df} 

The following lemma is a consequence of the definitions in \cite[Section 5.1]{hm11}.

\begin{lem}[\cite{hm11}]{\label{l:ComparingWeakTrainTracks}} Let $\varphi \in Out(F_r)$ be nongeometric fully irreducible and let $T, T'$ be weak train tracks for $\varphi$ with $\Lambda$-isometries $f_T \colon T \to T_+$ and $f_{T'} \colon T' \to T_+$. Then:
\begin{description}
\item[A] If $f_T$ and $f_{T'}$ are each injective on principal vertices, then they split equally.
\item[B] If $f_T$ is injective on principal vertices, then $f_{T'}$ splits at least as much as $f_{T}$.
\end{description}
\end{lem}

We will also use \cite[Proposition 5.4]{hm11}, which we record here as Proposition \ref{P:5.4}.

\begin{prop}[\cite{hm11}]\label{P:5.4} Let $\varphi \in Out(F_r)$ be nongeometric fully irreducible. Then for any train track representative $g \colon \Gamma \to \Gamma$ for $\varphi$ with associated $\Lambda$-isometry $g_{\infty} \colon \tilde{\Gamma} \to T_+$, there exists an $\varepsilon > 0$ so that, if $f \colon T \to T_+$ is any $\Lambda$-isometry, if $g_{\infty}$ splits at least as much as $f$, and if $Len(T) \leq \varepsilon$, then there exists a unique equivariant edge-isometry $h \colon \tilde{\Gamma} \to T$ such that $g_{\infty} = f \circ h$. Moreover, $h$ is a $\Lambda$-isometry. 
\end{prop}

\vskip10pt

\section{The Stable Axis Bundle}{\label{S:SAB}}

As mentioned above, the stable axis bundle was introduced in \cite{hm11} Subsection 6.5 as an object of interest and is studied here as a means to a more general proof of our main theorem. For this purpose we establish here rigorously properties previously believed true.

\begin{df}[Stable axis bundle]
Let $\psi \in Out(F_r)$ be ageometric fully irreducible. Then 
\begin{equation}\label{e:st}
ST(\psi) :=\{\Gamma \in TT(\psi) \mid \exists \text{ a fully stable train track representative } g \colon \Gamma \to \Gamma \text{ for } \psi \}.
\end{equation}

\noindent The \emph{stable axis bundle} is 
\begin{equation}\label{e:sab}
\mathcal{SA}_{\varphi} = \overline{\cup_{k=1}^{\infty} ST(\varphi^k)}.
\end{equation}
\end{df}

One can reformulate the stable axis bundle definition in terms of principal points: 

\begin{lem}{\label{L:SABReformulation}}
Let $g \colon \Gamma \to \Gamma$ be a rotationless train track representative of a positive power of an ageometric fully irreducible $\varphi \in Out(F_r)$. Then the following are equivalent: 
\begin{enumerate}
\item $\varphi$ is fully stable.
\item $\varphi$ has no Nielsen paths.
\item The associated map $g_{\infty} \colon \tilde{\Gamma} \to T_+$ is injective on the set of principal points.  
\end{enumerate}
\end{lem} 

\begin{proof}
That (1) implies (2) is simply Definition \ref{D:Ageometric}. (2) implies (1) by the definition of a stable train track representative.
We now show (2) implies (3). Assume $g$ is NP-free and that (3) does not hold, i.e. that there exist distinct points $\tilde{x},\tilde{y} \in \tilde{\Gamma}$ such that $g_{\infty}(\tilde{x}) = g_{\infty}(\tilde{y})$. Then Proposition \ref{P:C2.14} tells us that we can take a path $\sigma$ from $\tilde{x}$ to $\tilde{y}$, project to $\Gamma$, and know that some $g^k_{\#}(\sigma)$ will either be trivial or an NP. However, it cannot be trivial since the endpoints are distinct and periodic. Hence, some power of $\sigma$ must be an NP, contradicting our assumption.
We now prove that (3) implies (2). Suppose that $\Gamma$ has an NP $\sigma$. Let $\tilde{g}$ and $\tilde{\sigma}$ be lifts of $g$ and $\sigma$, respectively, so that $\tilde{g}$ preserves the endpoints $\tilde{x},\tilde{y}$ of $\tilde{\sigma}$. Then Proposition \ref{P:C2.14} implies that $g_{\infty}(\tilde{x}) = g_{\infty}(\tilde{y})$, i.e. that $g_{\infty}$ is not injective on the set of principal points. 
\end{proof}

\begin{df}[Stable weak train track]
Generalizing the lemma, one can define a weak train track $\Gamma$ to be \emph{stable} if there exists a $\Lambda$-isometry $\tilde{\Gamma} \to T_+$ which is injective on the set of principal points. We denote the set of stable weak train tracks for a given fully irreducible $\varphi$ by $SWTT(\varphi)$.
\end{df}

Proposition \ref{p:stts} will then imply that $\mathcal{SA}_{\varphi}$ and $SWTT(\varphi)$ are in fact the same set.

\begin{df}[Weak periodic Nielsen path]
A \emph{weak periodic Nielsen path} in a weak train track $T$ is a homotopically nontrivial path in $T$ whose endpoints are principal points with the same image in $T_+$. 
\end{df}

One can then also characterize stable weak train tracks by their lack of weak PNPs:

\begin{lem}
A weak train track $T$ is stable if and only if it has no weak PNPs.
\end{lem}

\begin{proof}
If $T$ is a stable weak train track, then it is injective on principal points. Hence, the endpoints of a weak PNP in $T$ would be the same, but this is impossible since $T$ is a tree. 

Suppose that $T$ is a weak train track with distinct principal points $v_1$, $v_2$ in $T$ having the same image in $T_+$. Since $T$, being a tree, is connected, there exists some path from $v_1$ to $v_2$ in $T$. This path would be a weak PNP. Hence, if $T$ is a weak train track, it is injective on principal points.
\end{proof}

\begin{lem}{\label{l:uniqueisometry}} Suppose $\varphi \in Out(F_r)$ is ageometric, fully irreducible with $i(\varphi) = \frac{3}{2}-r$. Then, for each $X \in \cup ST(\varphi^k)$, the $\Lambda$- isometry $\mathcal{I} \colon \tilde{X} \to T_{+}$ of the universal cover is unique.
\end{lem}

\begin{proof} Let $\varphi \in Out(F_r)$ be ageometric fully irreducible with attracting lamination $\Lambda$, i.e. $\Lambda = \Lambda_{\varphi}$. Consider a point in Outer Space viewed as a free, simplicial $F_r$-tree $T$. In \cite{hm11}, Handel and Mosher define an \emph{orientation} of $\tilde{\Lambda}$ in $T$ as an $F_r$-equivariant choice of orientation on each leaf of $\tilde{\Lambda}$ satisfying that the orientations of leaves $L_T$ and $L'_T$ agree on their intersection $L_T \cap L'_T$, provided the intersection contains a nontrivial interval. One can note that, by birecurrence, an orientation on $\tilde{\Lambda}$ is determined by an orientation of any of its leaves. One can also note that an orientation of $\tilde{\Lambda}$ in $T$ induces a well defined $F_r$-equivariant orientation on each edge of $T$, hence on its quotient graph $\Gamma=T/F_r$.

In \cite[Theorem 5.8]{hm11} they prove that the $\Lambda$-isometry $\mathcal{I} \colon T \mapsto T_+$ is unique if $\tilde{\Lambda}$ is nonorientable. Hence, we can assume $\tilde{\Lambda}$ is orientable and fix an orientation.

By a \emph{positive} gate we will mean a gate that the attracting lamination only exits (and never enters) the vertex through. On the other hand, a \emph{negative} gate will mean a gate that the attracting lamination only enters (and never exits) the vertex through. Each gate is either positive or negative. A direction in a positive gate will be called \emph{positive} and a direction in a negative gate will be called \emph{negative}. Notice that, for an orientable lamination realized in $\Gamma$, for each edge of $\Gamma$, the direction of the edge at one vertex is positive and the direction of the edge at the other vertex is negative.

Let $T$ specifically represent a point in $\cup ST(\varphi^k)$ and let $X=T/F_r$. Thus, there exists a fully stable train track representative $g\colon X \to X$ of some rotationless power $\varphi^R$. According to~\cite[Theorem 5.8c(iii)]{hm11}, $\mathcal{I} \colon T \to T_{+}$ is unique if and only if there exist vertices $v,w \in T$ such that $v$ has $\geq 2$ positive gates with respect to $\mathcal{I}$ and $w$ has $\geq 2$ negative gates with respect to $\mathcal{I}$. Thus, by symmetry, it suffices to prove it impossible for each vertex to have a unique positive gate.

For the sake of contradiction we assume each vertex has only one positive gate. We will use the fact that, in the absence of PNPs, the rotationless index is computed from the gates. In the situation where $i(\varphi) = \frac{3}{2}-r$, since $g$ is fully stable, ignoring vertices of valence 2, there is a single vertex $x$ of $X$ that contains precisely one nonperiodic direction and every direction of every other vertex is periodic. Notice, in particular, that each vertex other than $x$ (and having $\geq 3$ gates) has $\geq 3$ gates, so $\geq 2$ negative gates. Also, the total number of positive vertex directions in $X$ and of negative vertex directions in $X$ must be equal in order for them to correspond to a set of edge orientations.

We restrict our attention to vertices of valence $\geq 3$. We consider separately the cases where $X$ has only one such vertex and where $X$ has more than one such vertex. If $X$ has only one such vertex, then it would have to have $2r-2$ negative gates. Thus, it would have $\geq 2r-2$ negative directions and $\leq 2$ positive directions. For $r \geq 3$ this makes having an equal number of positive and negative directions impossible. So suppose $X$ has $k \geq 2$ such vertices. Then $X$ would have $\geq 2k$ negative gates by the previous paragraph. Thus, it would have $\geq 2k$ negative directions and at most $k+1$ positive directions. This is a contradiction, as above, unless $k=1$.
\qedhere
\end{proof}

Before proceeding with Proposition \ref{p:stts}, we recall from ~\cite[Proposition 5.5]{hm11} the following (the notation is described in Definition \ref{d:WhiteheadGraphs} and Definition \ref{d:Principal}).

\begin{prop}[~\cite{hm11}]{\label{p:Prop5.5}}
Let $T$ be a weak train track and $f \colon T \to T_+$ a $\Lambda$ - isometry. Then the following are necessary and sufficient conditions for the existence of a train track representative $g \colon \Gamma \to \Gamma$ of $\varphi^k$, for some $k$, such that $\Gamma = T /F_r$ and $f = g_{\infty}$ (i.e. for $\Gamma$ to be a train track for $\varphi^k$):
\begin{enumerate}
\item For every vertex $w$ of $T$ , $f(w) \in T_+$ is a periodic point.
\item For every vertex $y$ of $T$, if $f(y) = b$ is a branch point of $T_+$ then there exists a principal vertex $w$ of $T$ such that $Df (\mathcal{LW}(y;T)) \subset Df (SW (w; T )) \subset \mathcal{C}_b$.
\end{enumerate}
\end{prop}



\begin{prop}{\label{p:stts}} Suppose that $\varphi \in Out(F_r)$ is an ageometric, fully irreducible outer automorphism with $i(\varphi) = \frac{3}{2}-r$. Then
\begin{itemize}
\item[A.] $\mathcal{SA}_{\varphi}$ is the set of stable weak train tracks for $\varphi$
\item[B.] $\mathcal{SA}_{\varphi}$ is a nonempty, closed, $\varphi$-invariant subset of $\mathcal{A}_{\varphi}$.

\end{itemize}
\end{prop}

\begin{proof}
We first show, via applying Proposition \ref{p:Prop5.5}, that $SWTT(\varphi)$ is contained in $\overline{\cup ST(\varphi^k)}$. We cannot directly apply Proposition \ref{p:Prop5.5} to a given $T \in SWTT(\varphi)$ (with $\Lambda$-isometry $f_T \colon T \to T_+$) because $T$ may fail to satisfy the first of the necessary and sufficient conditions of each vertex being preperiodic. Hence, we approximate $T$ by performing the operation in the proof of ~\cite[Corollary 5.6]{hm11} of eliminating nonpreperiodic vertices $w$ one at a time via small folds of directions at $w$ having the same image. In fact, the folds can be chosen sufficiently small to avoid interaction with any principal vertices. In particular, for the $T'$ obtained, the injectivity of $f_{T'}$ on the set of principal vertices is unaffected, and the weak train tracks $f_{T'} \colon T' \to T_+$ obtained will still stable. Thus, we can apply Proposition \ref{p:Prop5.5} to approximate $T$ by stable weak train tracks provided that the $T'$ satisfy Proposition \ref{p:Prop5.5}(2), and hence are train tracks. That is, we need for each vertex $y$ of $T'$, such that $f(y)=b$ for some branch point of $T_+$, there exists a principal vertex $w$ of $T'$ such that 
\begin{equation}{\label{E:Containment}}
Df_{T'} (W (y; T' )) \subset Df_{T'} (SW (w; T' )) \subset \mathcal{C}_b.
\end{equation}
First notice that, because $f_{T'} \colon T' \to T_+$ is a stable weak train track, for each principal point $w \in T'$ and $b=f_{T'}(w)$ we have that $Df_{T'} (SW (w; T' )) = \mathcal{C}_b$. This follows from the fact that $$Df_{T'} (SW (w; T' )) \subset \mathcal{C}_b \subset \mathop{\cup}_{\text{principal } x} Df_{T'} (SW (x; T' ))$$ and that since $f_{T'}$ is injective on principal vertices, no leaf of 
$$\mathop{\cup}_{\text{principal } x} Df_{T'} (SW (x; T' )) - Df_{T'} (SW (w; T' ))$$ 
can be contained in $\mathcal{C}_b$. $Df_{T'} (W (y; T' )) \subset \mathcal{C}_b$ since it contains all leaves passing through $b$. Hence, 
(\ref{E:Containment}) holds and Proposition \ref{p:Prop5.5}(2) is satisfied. So the $T'$ are both stable weak train tracks and trains tracks and hence are stable train tracks by Lemma \ref{L:SABReformulation}. Since they approximate $T$, we have that $T \in \overline{\cup ST(\varphi^k)}$. Hence, $SWTT(\varphi) \subset \overline{\cup ST(\varphi^k)}$, as desired.


Since each stable train track is in $SWTT(\varphi)$, we are left to show for (A) that $SWTT(\varphi)$ is closed. In other words, we need that each $T \in \overline{SWTT(\varphi)}$ is in fact in $SWTT(\varphi)$, i.e. that the associated $\Lambda$-isometry $T \to T_+$ is injective on principal points. Notice that it also suffices to show this in $\widehat{CV_r}$, as the projection to $CV_r$ will then also be closed. Sometimes in what follows we will use the same notation for a tree and its projection.

Let $T$ be in the closure and $T_i$ a sequence of stable weak train tracks converging to $T$. Take a subsequence, if necessary, so that all $T_i$ are in the same open cell. Notice that, if $T$ is not in the open cell containing the $T_i$, then it is in a face of the cell. Let $f_i^{+}$ denote the $\Lambda$-isometry $T_i \to T_+$.



Let $X=\cup \{T_i\}\cup T$. Then $X$ is a compact subset of $\widehat{\mathcal{A}_{\varphi}}$. Hence, since the length function is continuous, there is an upper bound on the length of a $F_r$-tree in $X$. By Proposition \ref{P:6.2}, for each $F_r$-tree $R$ in $\widehat{\mathcal{A}_{\varphi}}$ and each integer $k$, we have $Len(\varphi^k(R)) = \lambda(\varphi)^{-k} Len(R)$. Thus, given any $\varepsilon$, there exists a $k_{\varepsilon}$ such that, for all $i$, we have $Len(\varphi^{k_{\varepsilon}}(T_i)) < \varepsilon$. Because applying $\varphi^{k_{\varepsilon}}$ does not change full stability (just acts as a change of marking) it is safe in what follows to replace $\varphi^{k_{\varepsilon}}(T_i)$ with $T_i$, to replace $\varphi^{k_{\varepsilon}}(T)$ with $T$, etc. 

By the arguments of the previous paragraph, we can assume that $Len(T_i), Len(T) < \varepsilon$ and apply Proposition \ref{P:5.4} as follows. Since the $\Lambda$-isometry for each $F_r$-tree in $\{T_i\}$ is injective on principal points, we can choose some $F_r$-tree $S$ in $\mathcal{A}_{\varphi} \cap (\cup ST(\varphi^k))$, with $\Lambda$-isometry $s \colon S \to T_+$ (also injective on principal points), such that $S$ satisfies: For each $T_i$, there exists a $\Lambda$-isometry $l_i\colon S \to T_i$ such that each arrow in the following diagram represents a $\Lambda$-isometry.

$$\xymatrix{S \ar[r]_{l_i} \ar@/^2pc/[rr]_{s} & T_i \ar[r]_{f_i^{+}} & T_{+} \\}$$

Since $T$ is in the closure of the cell containing the $T_i$, one can obtain $T$ from each $T_i$ via a quotient map $q_i \colon T_i \to T$, affine on edges. Let $m_i=q_i \circ l_i$.

$$\xymatrix{S \ar[r]_{l_i} \ar@/^2pc/[rr]_{m_i} & T_i \ar[r]_{q_i} & T \\}$$

For each $i$, let $G_i := T_i / F_r$ denote the quotient graph of $T_i$ and let $q_i' \colon G_i \to G$ be the induced quotient map. Since all $T_i$ are in the same open cell, there exists a family of marked homeomorphisms $g_{ji} \colon G_i \to G_j$, affine on each edge, so that $g_{kj} \circ g_{ji} = g_{ki}$. For each $G_i$ choose an indexing $\{e^{\alpha}_i\}$ of the edge set of $G_i$ so that $g_{ji}(e^{\alpha}_i) = e^{\alpha}_j$. In a well-defined manner we can give each edge $q_i'(e^{\alpha}_i)$ in $G$ the label $e^{\alpha}$. (Note that, for some $e^{\alpha}_i$, we have that $q_i'(e^{\alpha}_i)$ is just a point, but we are only concerned with cases where it is an edge.)

Since each $l_i$ is a $\Lambda$-isometry, and the lengths of the $\{e^{\alpha}_i\}$ converge to the lengths of the $e^{\alpha}$, the $m_i$ converge to a $\Lambda$-isometry $m\colon S \to T$. Then $g \circ m \colon S \to T_+$ is a composition of $\Lambda$ isometries, hence is a $\Lambda$-isometry. By Lemma \ref{l:uniqueisometry}, $s=g \circ m$.

$$\xymatrix{S \ar[r]_{m} \ar@/^2pc/[rr]_{s} & T \ar[r]_{g} & T_{+} \\}$$

Recall that, by definition, principal points of an $F_r$-tree $R$ in $\mathcal{A}_{\varphi}$ map via the $\Lambda$-isometry $R \to T_{+}$ to branch points in $T_{+}$. What we need to show is that the set $\mathscr{S}(T)$ of principal points of $T$ is mapped injectively via $g$ into the set of branch points of $T_{+}$, and $T \in SWTT(\varphi)$, as desired. So, in particular, we need to show that, given principal points $v,w \in T$ with $g(v) = g(w)$, we have $v = w$. 

Let $v,w \in T$ be principal points such that $g(v) = g(w)=b$ and let $l$ and $l'$ be leaves of $\mathcal{C}_b$ such that $v$ is the basepoint of the realization $l_T$ of $l$ in $T$ and $w$ is the basepoint of the realization $l'_T$ of $l'$ in $T$. Let $b_S$ be the basepoint of the realization $l_S$ of $l$ in $S$ and let $b'_S$ be the basepoint of the realization $l'_S$ of $l'$ in $S$. Then $b_S = b'_S$, as $s$ is injective on principal points and $s(b_S) = b = s(b'_S)$. Now, since $m$ is a $\Lambda$-isometry, $m$ maps $l_S$ isometrically onto $l$ and $l'_S$ isometrically onto $l'$. Also, $m(b_S) = v$ since $m$ maps $l_S$ isometrically onto $l_T$ and $v$ is the only point on $l_T$ mapped to $b$. Similarly, $m(b'_S) = w$. Thus, $v = w$, as desired. And the proof that $SWTT(\varphi)$, and hence $\mathcal{SA}_{\varphi}$, is closed is complete.

$SWTT(\varphi)$ is nonempty since $\varphi$ has a rotationless power with a stable representative (obtained, for example, via the stabilization algorithm of \cite{bh92}).
Since $\mathcal{SA}_{\varphi}$ is a closure, hence closed, we are left to show that $SWTT(\varphi)$ is $\varphi$-invariant. That $\cup ST(\varphi^k)$ is invariant follows from the fact that changing the marking does not change stability. $\mathcal{SA}_{\varphi}$ is then invariant since the action is continuous.
\qedhere
\end{proof}

\vskip10pt

\section{The Proof}{\label{S:Proofs}}

According to \cite{hm11}, $\mathcal{A}_{\varphi}$ is proper homotopy equivalent to a line. In fact, there are distinguished lines from which it is possible to get from to any point in $\mathcal{A}_{\varphi}$ by a combination of folding and of collapsing PNPs. These distinguished lines are periodic fold lines for representatives with the maximum number of PNPs possible. The strategy of our proof is to show when $\mathcal{A}_{\varphi}$ consists of only a single such fold line.

We start with a sequence of two lemmas immediately revealing the necessity for the rotationless index to be $\frac{3}{2}-r$.

\begin{lem}{\label{l:ec}} Suppose $g \colon \Gamma \to \Gamma$ is a fully stable train track representative of a rotationless power $\varphi^R$ of a fully irreducible $\varphi \in Out(F_r)$. Then:
\begin{description}
\item[A] If the rotationless index satisfies $i(\varphi) = \frac{3}{2}-r$, then $g$ has a unique illegal turn.
\item[B] If the rotationless index satisfies $i(\varphi) > \frac{3}{2}-r$, then by following from $\Gamma$ a Stallings fold decomposition of $g$, one reaches a point $\Gamma' \in CV_r$ such that there exists a fully stable train track representative $g' \colon \Gamma' \to \Gamma'$ of some power $\varphi^R$ with more than one illegal turn.
\end{description}
\end{lem}

\begin{proof}
By \cite{gjll}, $\varphi$ is ageometric if and only if $i(\varphi) \geq \frac{3}{2}-r$. We thus can assume each stable representative is PNP-free. 
 
For simplicity, in what follows, we denote 
$$GI(g) = \sum\limits_{\text{vertices v}}(1-\frac{\#\text{(gates of v)}}{2}).$$ 
Recall also (from (Definition \ref{d:IndexLists})) that the rotationless index $i(\varphi)$ is the same sum as $GI(g)$ but with the terms indexed by nonperiodic vertices with $\geq 3$ gates removed. Hence, in particular, $GI(g) \leq i(\varphi)$.


The nonzero terms appearing in the sum for $GI(g)$, but not for $i(\varphi)$, are precisely those indexed by preperiodic vertices v with $\geq 3$ gates.
  
$$1-r = \chi(\Gamma) = \# Vertices - \#Edges = \# Vertices - \frac{1}{2}\sum\limits_{\text{gates D}}card(D)$$ 
$$= GI(g) + \frac{1}{2}\sum\limits_{\text{gates D}}(1-card(D))$$
Thus,
\begin{equation}{\label{E:EulerEquation}}
(1-r) - GI(g) = \frac{1}{2}\sum\limits_{\text{gates D}}(1-card(D)).
\end{equation}

Each term on the right-hand side of (\ref{E:EulerEquation}) is nonpositive and is zero if and only if the gate $D$ consists of a single direction. Thus, $GI(g)=1-r$ if and only if each gate has cardinality one. This is true if and only if there are no illegal turns, which is impossible as $\varphi$ is fully irreducible, hence of infinite order. Thus, $GI(g)>1-r$ and, since $GI(g)$ can only take $1/2$-integer values, $GI(g) \geq \frac{3}{2}-r$. Hence, under the assumption that $i(\varphi) = \frac{3}{2}-r$, we have
$$\frac{3}{2}-r \leq GI(g) \leq i(\varphi) = \frac{3}{2}-r.$$ 
Now, $GI(g) = \frac{3}{2}-r$ precisely when each gate has cardinality 1, except a single gate of cardinality 2. Equivalently, $GI(g) = \frac{3}{2}-r$ if and only if $g$ has a unique illegal turn, proving A.

We suppose $i(\varphi) > \frac{3}{2} - r$ and prove B. The first observation we use is that $GI(g) \leq i(\varphi)$, with equality if and only if each vertex with $\geq 3$ gates is fixed. The second observation is that the following three statements are equivalent:

\noindent (i) $GI(g) > \frac{3}{2}-r$.

\noindent (ii) Either there exist gates $D_1$, $D_2$ with card($D_1$), card($D_2$) $>1$ or there exists some gate $D$ with card($D$) $>2$.

\noindent (iii) $g$ has $\geq 2$ illegal turns.

Let $g \colon \Gamma \to \Gamma$ be any rotationless train track representative of a $\varphi^R$. If $g$ has $\geq 2$ illegal turns, then B is proved. So we suppose $g$ has only one illegal turn. As above, $GI(g) = \frac{3}{2} - r$. So $i(\varphi) > GI(g)$ and $\Gamma$ must have a nonperiodic vertex with $\geq 3$ gates. Any such vertex $v$ is preperiodic, i.e. there exists a fixed vertex $w$ and $j \geq 1$ such that $g^j(v) = g^j(w) = w$.

By iterating a Stallings fold decomposition of $g \colon \Gamma \to \Gamma$ we obtain periodic fold line as in Subsection \ref{SSS:PeriodicFoldLines} (see Display \ref{E:PeriodicFoldLines}, in particular).

Let $k$ be such that $g_{k+1} \circ \dots \circ g_1(v) = g_{k+1} \circ \dots \circ g_1(w)$, while $g_{k} \circ \dots \circ g_1(v) \neq g_{k} \circ \dots \circ g_1(w)$. For each $i$ let $v_i$ denote $g_{i} \circ \dots \circ g_1(v)$, let $w_i$ denote $g_{i} \circ \dots \circ g_1(w)$, and let $f_i=g_{i} \circ \dots \circ g_1 \circ g_K \circ \dots \circ g_{i+1} \colon \Gamma_i \to \Gamma_i$. Notice that the fold map $g_{k} \colon \Gamma_k \to \Gamma_{k+1}$ conjugates $f_k$ to $f_{k+1}$, i.e. $g_{k+1} f_k = f_{k+1} g_{k+1}$. Also, the $g_k$ bijectively map the periodic directions of the $w_i$ because periodic directions cannot be identified (as they are in distinct gates) and $g$ bijectively maps the periodic directions.

In $\Gamma_k$, the vertex $v_k$ is preperiodic (and not periodic) and has $\geq 3$ gates, while the vertex $w_k$ is fixed with $\geq 3$ gates. To identify $v_k$ and $w_k$, the fold $g_{k+1} \colon \Gamma_k \to \Gamma_{k+1}$ must be an improper full fold, which fully folds two oriented edges $E$, $E'$ having the same initial vertex and having terminal vertices $v_k$, $w_k$ respectively. We have $g_{k+1}(E)=g_{k+1}(E')=E''$ with terminal vertex $g_{k+1}(v_k) = g_{k+1}(w_k) = w_{k+1}$. At $v_k$ there are $\geq 2$ directions, namely $d$ and $d'$, not in the same gate as each other nor in the same gate as the terminal direction of $E$. It follows that $g_{k+1}(d)$ and $g_{k+1}(d')$ are two directions at $w_{k+1}$ not in the same gate as each other, and they are not periodic directions. Therefore, the two gates at $w_{k+1}$ containing $g_{k+1}(d)$ and $g_{k+1}(d')$ must each contain $\geq 1$ other direction, namely some periodic direction. Therefore, the train track representative $f_{k+1}$ has $\geq 2$ gates each of cardinality $\geq 2$. This proves $g' = f_{k+1} \colon \Gamma_{k+1} \to \Gamma_{k+1}$ is the desired fully stable train track representative with more than one illegal turn.
\qedhere
\end{proof}

\begin{com}{\label{c:c1}} Notice that item \textbf{B} does \emph{not} assert that $\Gamma$ itself has more than one illegal turn. This motivates a question which may shed light on whether local dimension is constant on the axis bundle: does there exist a fully irreducible $\varphi \in Out(F_r)$ and two fully stable train track representatives $g_1 \colon \Gamma_1 \to \Gamma_1$ and $g_2 \colon \Gamma_2 \to \Gamma_2$, such that $\Gamma_1$ has only one illegal turn and $\Gamma_2$ has more than one?
\end{com}

\begin{lem}{\label{l:distinct}} Suppose that $\varphi \in Out(F_r)$ is ageometric and fully irreducible. If $i(\varphi) > \frac{3}{2}-r$, then, for each $k \geq 1$, each point in $\cup ST(\varphi^k)$ is contained in at least two distinct periodic fold lines.
\end{lem}

\begin{proof} Suppose $\varphi \in Out(F_r)$ is ageometric and fully irreducible with rotationless index $i(\varphi) > \frac{3}{2}-r$. Then each $\varphi^k$ is also ageometric and fully irreducible with the same rotationless index. Let $\Gamma$ be a point of $\cup ST(\varphi^k)$ and $f \colon \Gamma \to \Gamma$ a fully stable rotationless train track representative of some $\varphi^k$. By Lemma \ref{l:ec}, we can fold from $\Gamma$ to an $X$ on which there exists a fully stable representative $g$ of $\varphi^k$ with more than one illegal turn, and that is in fact conjugate to $f$ by precisely the folds taken to move from $\Gamma$ to $X$. 

For any illegal turn $\{d, d'\}$ in $X$ consider its forward orbit $\{Dg^i(d),Dg^i(d')\}_{i \geq 0}$. There is a minimal value of $i \geq 1$ for which the turn $\{Dg^i(d), Dg^i(d')\}$ is degenerate. It follows that there exists an illegal turn $\{d,d'\}$ in $X$ whose immediate forward image $\{Dg(d), Dg(d')\}$ is degenerate. If there are two such illegal turns then, by choosing each of those two turns (respectively) for the first fold, we get two different Stallings fold factorizations of $g$, hence two different fold lines passing through $X$. Now suppose there is only one illegal turn $T_1 = \{d, d'\}$ such that $\{Dg(d), Dg(d')\}$ is degenerate. Since some other illegal turn exists, there exists a distinct illegal turn $T' = \{d'', d'''\} \neq T_1$ which is mapped to $T_1$ by $Dg$. For the first fold line passing through $X$, one can use a Stallings fold decomposition for $g$ starting with a maximal fold of $T_1$. One obtains a second fold line as follows. Start by folding the two initial segments corresponding to $T_1$ some nonmaximal amount, producing a fold segment from $X$ to some $X'$. Let $s$ denote the direction in $X'$ of the folded segment and $S$ the segment in the direction of $s$ that was obtained by the initial fold of $T_1$. Now we must consider separately the cases where (1) $T_1$ and $T'$ share no common direction and (2) where $T_1$ and $T'$ share a common direction. If (1) holds, we continue by folding the initial segments of $T'$ mapping to $S$, return to finish maximally folding $T_1$, then continue with any Stallings fold decomposition from there. Now assume (2) holds and $T_1$ and $T'$ share a common direction, say $d=d'''$. Let $d''$ denote the image in $X'$ of the direction of $X$ of the same name. Next fold two initial segments corresponding to the turn $\{s, d''\}$. Again one can continue and obtain a different Stallings fold decomposition, hence a different periodic line passing through $X$.

To obtain two distinct periodic fold lines passing through $\Gamma$, we create periodic fold lines for $\varphi^2$ by folding along the Stallings fold decomposition for $f$ from $\Gamma$ to $X$, following one of the distinct fold lines from $X$ to $X$, and then finishing from $X$ to $\Gamma$ the Stallings fold decomposition for $f$.
\qedhere
\end{proof}

Another obstacle to an axis bundle having only one axis is the possible existence of multiple affine train track representative on the same point of Outer Space. The following lemma restricts when this can occur (by Proposition \ref{P:C2.14} each train track representative would induce a distinct $\Lambda$-isometry of the universal cover).


\begin{lem}{\label{l:uniqueisometry2}} Suppose that $\varphi \in Out(F_r)$ is ageometric fully irreducible with $i(\varphi) = \frac{3}{2}-r$. Suppose that $X, Y \in \cup ST(\varphi^k)$ and that there exists a $\Lambda$-isometry $\mathcal{I} \colon \tilde{X} \to \tilde{Y}$. Then $\mathcal{I}$ is unique.
\end{lem}

\begin{proof} Suppose there were two distinct $\Lambda$-isometries $\mathcal{I}_1,\mathcal{I}_2 \colon \tilde{X} \to \tilde{Y}$. By Lemma \ref{l:uniqueisometry}, there is a unique $\Lambda$- isometry $\mathcal{I}_X \colon \tilde{X} \to T_+$ and unique $\Lambda$- isometry $\mathcal{I}_Y \colon \tilde{Y} \to T_+$. Let $\mathcal{I}_{X,1} = \mathcal{I}_Y \circ \mathcal{I}_1 \colon \tilde{X} \to T_+$ and $\mathcal{I}_{X,2} = \mathcal{I}_Y \circ \mathcal{I}_2 \colon \tilde{X} \to T_+$. Let $x \in \tilde{X}$ be such that $\mathcal{I}_1(x) \neq \mathcal{I}_2(x)$. Let $L$ be a leaf of $\Lambda$ realized in $\tilde{X}$ and passing through $x$. Since $\mathcal{I}_1$ and $\mathcal{I}_2$ are $\Lambda$- isometries, $\mathcal{I}_1(L) = \mathcal{I}_2(L)$. Thus, for each $y \in L$, we have that $\mathcal{I}_1(y)$ is a shift of $\mathcal{I}_2(y)$ along $L$ by the same distance as $\mathcal{I}_1(x)$ is shifted from $\mathcal{I}_2(y)$. Since $\mathcal{I}_Y$ is a $\Lambda$- isometry, this also holds for $\mathcal{I}_{X,1}$ and $\mathcal{I}_{X,2}$. This contradicts the uniqueness of $\mathcal{I}_X$. 
\qedhere
\end{proof}

\begin{lem}{\label{l:sisa}} Suppose that $\varphi \in Out(F_r)$ is ageometric fully irreducible such that $i(\varphi) = \frac{3}{2}-r$. Then $\mathcal{IW}(\varphi)$ has no cut vertices if and only if $\mathcal{SA}_{\varphi}=\mathcal{A}_{\varphi}$.
\end{lem}

\begin{proof}
Suppose no component of $\mathcal{IW}(\varphi)$ has a cut vertex. Then, by ~\cite[Lemma 3.1]{hm11}, no train track representative of $\varphi$ has a PNP. That is, every representative of $\varphi$ is fully stable. Since $\mathcal{IW}(\varphi) = \mathcal{IW}(\varphi^k)$ for all $k \geq 1$, the same can be said for each $\varphi^k$ with $k \geq 1$. Thus, $TT(\varphi^k) = ST(\varphi^k)$ for each $k \geq 1$ and $\mathcal{SA}_{\varphi}=\mathcal{A}_{\varphi}$, as desired.

Suppose $\mathcal{IW}(\varphi)$ has a component with a cut vertex. We claim that ~\cite[Lemma 4.3]{hm11} implies that some power  $\varphi^k$ has a train track representative (we can consider to be affine) $g: \Gamma \to \Gamma$ with a PNP. In order to apply ~\cite[Lemma 4.3]{hm11}, we need to show that there exists a PNP-free train track representative $g' \colon \Gamma' \to \Gamma'$ of this sufficiently high power of $\varphi$ satisfying (1)-(3) of the lemma. We take a high enough power $k$ so that the image of every nonperiodic vertex is a fixed vertex and, further, the image of every nonperiodic direction is a fixed direction. We write this representative $g'' \colon \Gamma'' \to \Gamma''$. Since $g''$ has no PNP and the ideal Whitehead graph has a component with a cut vertex, there exists a principal vertex $w$ of $\Gamma''$ such that $\mathcal{SW}(g'',w)$ has a cut vertex. At the vertex $w$, we fold each gate at $w$ a small amount so that each direction at $w$ is fixed by our new train track representative $g' \colon \Gamma' \to \Gamma'$. Notice that now $\mathcal{SW}(g'',w) \cong \mathcal{SW}(g',w) \cong \mathcal{LW}(g',w)$, so that $\mathcal{LW}(g',w)$ has a cut vertex, which we will call $x$. 
We therefore have a decomposition into nontrivial subgraphs $\mathcal{LW}(g',w) = X_1 \cup X_2$ such that $X_1 \cap X_2 = \{x\}$, which verifies (1).
Since $g'$ is obtained from $g''$ by folding at illegal turns, $g'$ also has no PNPs. For (2), we need that if $g'(v)=w$ for some vertex $v \neq w$, then $Dg(\mathcal{LW}(g',v))$ is contained in either $X_1$ or $X_2$. By the proof of Lemma \ref{l:ec}, in order to have $i(\varphi)=\frac{3}{2}-r$, we need that $GI=\frac{3}{2}-r$, which can only occur when each preperiodic vertex either has valence 2 or has valence 3 and the unique illegal turn. If $v$ has valence 2, then $\mathcal{LW}(g',v)$ is a single edge, so its image is in the component. Also in the case where $v$ has valence 3, and the unique illegal turn, the image of $\mathcal{LW}(g',v)$ is a single edge. In either case, the single edge must lie in either $X_1$ or $X_2$, which completes the verification of (2). Since $\mathcal{SW}(g',w) = \mathcal{LW}(g',w)$ and, in particular, all directions at $w$ are fixed by $Dg'$, we have that $Dg'$ is the identity on $\mathcal{LW}(g',w)$, so that (3) hold. Since $g'(w')=w'$, we then have from ~\cite[Lemma 4.3]{hm11} the desired train track representative $g$ of $\varphi^k$ with a PNP.

This $\Gamma$ that $g$ is a train track map on, together with its marking and metric, gives a point in $\mathcal{A}_{\varphi}$. By Lemma \ref{l:uniqueisometry}, there cannot additionally be a stable train track representative on that point. Thus, the point is not in $\cup ST(\varphi^k)$. By Proposition \ref{p:stts}, $\mathcal{SA}_{\varphi}$ is actually the set of stable weak train tracks. And, hence, since a stable weak train track cannot be induced by a train track map unless it is a stable train track map (see Lemma \ref{L:SABReformulation}), $g$ cannot be in $\mathcal{SA}_{\varphi}$.
\qedhere
\end{proof}

\begin{thm}{\label{t:stablebundle}} Suppose $\varphi \in Out(F_r)$ is ageometric fully irreducible. Then the stable axis bundle $\mathcal{SA}_{\varphi}$ is a unique axis if and only if the rotationless index satisfies $i(\varphi) = \frac{3}{2}-r$. In that case it is a unique periodic fold line.
\end{thm}

\begin{proof}
Notice that $\mathcal{SA}_{\varphi}$ would have to contain the entirety of the periodic fold line for each affine train track representative on each element of $\cup ST(\varphi^k)$. Suppose that $\mathcal{SA}_{\varphi}$ contained more than one fold line. Since $\cup ST(\varphi^k)$ is dense in $\mathcal{SA}_{\varphi}$, then $\mathcal{SA}_{\varphi}$ would contain stable train tracks on distinct periodic fold lines. Hence, $\mathcal{SA}_{\varphi}$ would contain distinct periodic fold lines.

First suppose $i(\varphi) > \frac{3}{2}-r$. Then each point in $\cup ST(\varphi^k)$ is contained in at least two distinct periodic fold lines by Lemma \ref{l:distinct}. We thus suppose instead that $i(\varphi) = \frac{3}{2}-r$. It suffices to prove that $\mathcal{SA}_{\varphi}$ contains a unique periodic fold line.

Before proceeding with the proof, we remark that, for each $T \in \cup ST(\varphi^k)$, there is only one way to fold from $T$ to $T_+$. This is because there is only one $\Lambda$- isometry from $T$ to $T_+$ (by Lemma \ref{l:uniqueisometry}) and only one illegal turn (by Lemma \ref{l:ec}). 
 
We remark further that, for a given $T \in SWTT(\varphi)$, the realization of $\mathcal{IW}(\varphi)$ in $T$ is the disjoint union of the stable Whitehead graphs $\mathcal{SW}(w;T)$ for the principal points $w$ of $T$. Thus, for any $T,T' \in SWTT(\varphi)$ with respective $\Lambda$- isometries $i_T \colon T \to T_+$ and $i_{T'} \colon T' \to T_+$, by Lemma \ref{l:ComparingWeakTrainTracks}, we have that $i_T$ and $i_{T'}$ both split minimally (and, in particular, split as much as each other). This allows us to apply Proposition \ref{P:5.4}. That is, for any $T \in SWTT(\varphi)$, there exists an $\varepsilon>0$ such that, for any $T' \in SWTT(\varphi)$ with $Len(T') \leq \varepsilon$, the $\Lambda$- isometry $i_T \colon T \to T_+$ factors (uniquely) as a $\Lambda$- isometry $i_T \colon T \to T'$ followed by the $\Lambda$- isometry $i_{T'} \colon T' \to T_+$. Since there can only be one fold line from $T'$ to $T_+$, this implies $T'$ lies on the unique fold line from $T$ to $T_+$.

For the sake of contradiction suppose $\mathcal{SA}_{\varphi}$ contained two distinct periodic fold lines $L$ and $L'$. 
Choose $T \in \cup ST(\varphi^k)$ on $L$, choose $\varepsilon$ as in Proposition \ref{P:5.4}, and then choose $T' \in \cup ST(\varphi^k)$ on $L'$ so that $Len(T') < \varepsilon$. 
Since the lines converge to $T_+$ (and $\cup ST(\varphi^k)$ is dense in $\mathcal{SA}_{\varphi}$), this is possible. Then, by Proposition \ref{P:5.4}, as explained in the previous paragraph, $T$ and $T'$ must be on a common fold line. Without generality loss assume that, on the fold line, the parameter of $T$ is less than the parameter of $T'$. Thus, for the parameter $t_0$ for $T'$, we have $L(t)=L'(t)$ for all $t \geq t_0$.

We claim that there cannot exist two distinct periodic fold lines reparametrizable so that, for some $t_0$, $L(t)=L'(t)$ for all $t \geq t_0$. Suppose that $L$ is a periodic fold line for $\varphi^k$ and $L'$ is a periodic fold line for $\varphi^{k'}$. Letting $m=kk'$, it follows that $L,L'$ are both periodic fold lines $\varphi^m$. This means that $\varphi^{m\alpha}(L(t))=L(t+m\alpha\log(\lambda))$ and $\varphi^{m\alpha}(L'(t))=L'(t+m\alpha\log(\lambda))$ for all $t \in \RR$ and $\alpha \in \ZZ$. Given $t \in \RR$, this gives us $L(t)=L(t')$ by using a suitable choice of $t' \geq t_0$ and $\alpha \in \ZZ$. 

\qedhere
\end{proof}

\begin{thm}{\label{t:uniqueaxis}} The axis bundle of an ageometric, fully irreducible $\varphi \in Out(F_r)$ is a unique axis precisely if both of the following two conditions hold:
~\\
\vspace{-5mm}
\begin{enumerate}
\item the rotationless index satisfies $i(\varphi) = \frac{3}{2}-r$ and 
\item no component of the ideal Whitehead graph $\mathcal{IW}(\varphi)$ has a cut vertex.
\end{enumerate}
\end{thm}

\begin{proof}
If both conditions hold, $\mathcal{A}_{\varphi}$ is a unique axis by Lemma \ref{l:sisa} and Theorem \ref{t:stablebundle}.

Suppose $\varphi \in Out(F_r)$ is ageometric fully irreducible and that $\mathcal{A}_{\varphi}$ is a unique axis. Since $\varphi$ is ageometric, by \cite{gjll}, $i(\varphi) \geq \frac{3}{2}-r$. If $i(\varphi) > \frac{3}{2}-r$, then each point in $\cup ST(\varphi^k)$ is contained in at least two distinct periodic fold lines by Lemma \ref{l:distinct}. Since each $\varphi$ has a stable train track representative, $\cup ST(\varphi^k)$ is nonempty. $\mathcal{A}_{\varphi}$ would contain multiple fold lines. So $i(\varphi) = \frac{3}{2}-r$.

The second condition now follows by Lemma \ref{l:sisa}.
\qedhere
\end{proof}

\vskip10pt

\section{Final Remarks}{\label{S:Applications}}

With the Coulbois computer package, the full irreducibility criterion of Pfaff \cite{p12c}, and full irreducibility decidability algorithm of Kapovich \cite{k12}, it is becoming increasingly easy to check that an outer automorphism satisfies the conditions to have a single-axis axis bundle. Once one determines that it does, they can compute a train track representative, using the Coulbois computer package. A Stallings fold decomposition of the representative then gives the periodic line, which is the entire axis bundle.

\vskip10pt

\noindent \textbf{Ideal decomposition diagrams}

\vskip1pt

Let $\varphi \in Out(F_r)$ be ageometric fully irreducible. Suppose that, in addition to the conditions of Theorem \ref{t:uniqueaxis}, the ideal Whitehead graph is connected. Then $\mathcal{A}_{\varphi}$ is still a single periodic fold line. By \cite{p12b}, a power of $\varphi$ has a representative with a Stallings fold decomposition that is a sequence of proper full folds of roses. Also by \cite{p12b} we know that this decomposition has a realization as a loop in the ``ideal decomposition diagram'' $\mathcal{ID}(\mathcal{IW}(\varphi))$. In other words, the single axis in $\mathcal{A}_{\varphi}$ can be viewed as a repeated gradual folding of a loop in $\mathcal{ID}(\mathcal{IW}(\varphi))$. And this is true for any fully irreducible $\psi$ with $\mathcal{IW}(\psi) \cong \mathcal{IW}(\varphi)$.

\vskip10pt

\noindent \textbf{The conjugacy problem}

\vskip1pt

One can observe that $\mathcal{A}_{\varphi}$ and $\mathcal{A}_{\psi}$ differ by the action of $Out(F_r)$ on $CV_r$ if and only if there exist integers $k,l \geq 1$ such that $\varphi^k$ and $\psi^l$ are conjugate in $Out(F_r)$. Thus, given two outer automorphisms $\varphi$ and $\psi$ that one has checked, as above, satisfy the conditions for a single axis, one can construct the axis of each to determine if some $\varphi^k$ and $\psi^l$ are actually conjugate. In fact, $\varphi^k$ and $\psi^l$ are conjugate if they give the same bi-infinite path in the \cite{p12c} automata. 

\vskip10pt

\noindent \textbf{Determining all train track representatives for an outer automorphism}

\vskip1pt

In general, it is difficult to identify the set of train track representatives for a given outer automorphism. However, under the conditions of Theorem \ref{t:uniqueaxis}, once one has the axis bundle, all train track representatives have the same periodic fold line, namely the single axis of the axis bundle.

\vskip50pt

\bibliographystyle{amsalpha}
\bibliography{PaperReferences}

\end{document}